\documentclass[onecolumn]{article}
\usepackage{amssymb}
\usepackage{amsmath}
\usepackage{chicago}

\newtheorem{theorem}{Theorem}
\newtheorem{corollary}{Corollary}
\newtheorem{definition}{Definition}
\newtheorem{lemma}{Lemma}
\newtheorem{proposition}{Proposition}

\input{tcilatex}
\begin{document}

\title{The norm of products of free random variables}
\author{Vladislav Kargin \thanks{%
Courant Institute of Mathematical Sciences; 109-20 71st Road, Apt. 4A,
Forest Hills NY 11375; kargin@cims.nyu.edu}}
\date{}
\maketitle

\begin{abstract}
  Let $X_i$ denote free identically-distributed random variables. This paper
investigates how the norm of products $\Pi _n=X_1 X_2 ... X_n$ behaves as
$n$ approaches infinity. In addition, for positive $X_i$ it studies the
asymptotic behavior of the norm of $Y_n=X_1 \circ X_2 \circ ...\circ
X_n$, where $\circ$ denotes the symmetric product of two positive operators:
$A \circ B=:A^{1/2}BA^{1/2}$.

  It is proved that if the expectation of $X_i$ is 1, then the norm of the symmetric product $Y_{n}$ is
between $c_1 n^{1/2}$ \ and $c_2 n$ for certain constant $c_1$ and $c_2$. That is, the growth in the norm is at most linear.

 For the norm of the usual product $Pi_n$, it is proved that the limit of
$n^{-1}\log Norm(Pi_n)$ exists and equals $\log \sqrt{E\left( X_i^{\ast }X_{i}\right)}.$ In other words, the growth in the norm of the product is exponential and the rate equals the logarithm of the Hilbert-Schmidt norm of operator X.

Finally, if $\pi $ is a cyclic representation of the algebra generated by $X_i$, and if $\xi$ is a cyclic vector, then $n^{-1}\log Norm(\pi \left( \Pi _{n}\right) \xi)=\log \sqrt{E\left( X_{i}^{\ast }X_{i}\right) }$ for all $n.$
In other words, the growth in the length of the cyclic vector is exponential and the rate coincides with the rate in the growth of the norm of the product.
  
These results are significantly different from analogous results for
commuting random variables and generalize results for random matrices derived by Kesten and Furstenberg.

\end{abstract}

\section{Introduction}

Suppose $X_{1},$ $X_{2},$ ... $,$ $X_{n}$ are identically-distributed free
random variables. These variables are infinite-dimensional linear operators
but the reader may find it convenient to think of them as very large random
matrices. The first question we will address in this paper is how the norm
of $\Pi _{n}=X_{1}X_{2}...X_{n}$ behaves. If $X_{i}$ are all positive, then
it is natural to look also at the symmetric product operation $\circ $
defined as follows: $X_{1}\circ X_{2}=X_{1}^{1/2}X_{2}X_{1}^{1/2}.$ The
benefit is that unlike the usual operator product, this operation maps the
set of positive variables to itself. For this operation we can ask how the
norm of symmetric products $Y_{n}=X_{1}\circ X_{2}\circ ...\circ X_{n}$
behaves.\footnote{%
The operation $\circ $ is neither commutative, nor associative. By
convention we multiply starting on the right, so, for example, $X_{1}\circ
X_{2}\circ X_{3}\circ X_{4}=X_{1}\circ \left( X_{2}\circ \left( X_{3}\circ
X_{4}\right) \right) .$ However, this convention is not important for the
question that we ask. First, it is easy to check that $X_{1}\circ X_{2}$ has
the same spectral distribution and therefore the same norm as $X_{2}\circ
X_{1}.$ Second, if $X_{1},$ $X_{2},$ and $X_{3}$ are free, then the spectral
distribution of $\left( X_{1}\circ X_{2}\right) \circ X_{3}$ is the same as
the spectral distribution of $X_{1}\circ \left( X_{2}\circ X_{3}\right) ,$
and therefore these two products have the same norm. In brief, if $X_{i}$
are free, then the norm of $X_{1}\circ X_{2}\circ ...\circ X_{n}$ does not
depend on the order in which $X_{i}$ are multiplied by the operation $\circ
. $}

Products of random matrices and their asymptotic behavior were originally
studied by \shortciteN{bellman54}. One of the decisive steps was made by %
\shortciteN{furstenberg_kesten60}, who investigated a matrix-valued
stationary stochastic process $X_{1},$ $...$ $,$ $X_{n},$ $...$ , and proved
that the limit of $n^{-1}E\left( \log \left\Vert X_{1}...X_{n}\right\Vert
\right) $ exists (but might equal $\pm \infty $) and that under certain
assumptions $n^{-1}\log \left\Vert X_{1}...X_{n}\right\Vert $ converges to
this limit almost surely. Essentially, the only facts that are used in the
proof of this result are the ergodic theorem, the norm inequality $%
\left\Vert X_{1}X_{2}\right\Vert \leq \left\Vert X_{1}\right\Vert \left\Vert
X_{2}\right\Vert $ and the fact that the unit sphere is compact in
finite-dimensional spaces. It is the lack of compactness of the unit sphere
in the infinite-dimensional space that makes generalizations to
infinite-dimensional operators non-trivial (see \shortciteN{ruelle82} for a
generalization in the case of compact operators). More work on
non-commutative products was done by \shortciteN{furstenberg63}, %
\shortciteN{oseledec68}, \shortciteN{kingman73}, and others. The results are
often called multiplicative ergodic theorems and they find many applications
in mathematical physics. For example, see \shortciteN{ruelle84}.

In this paper, we study products of free random variables. These variables
are (non-compact) infinite-dimensional operators which can be thought of as
a limiting case of large independent random matrices.

Suppose that $X_{i}$ are free, identically-distributed, self-adjoint, and
positive. Suppose also $E\left( X_{i}\right) =1.$ Then we show that the norm
of $Y_{n}=X_{1}\circ X_{2}\circ ...\circ X_{n}$ grows no faster than a
linear function of $n.$ Precisely, we find that 
\begin{equation*}
\lim \sup_{n\rightarrow \infty }n^{-1}\left\Vert Y_{n}\right\Vert \leq
c_{1}\left\Vert X_{i}\right\Vert .
\end{equation*}%
We are also able to show that if $X_{i}$ is not concentrated at $1,$ then 
\begin{equation*}
\lim \inf_{n\rightarrow \infty }n^{-1/2}\left\Vert Y_{n}\right\Vert \geq
c_{2}>0.
\end{equation*}

For the usual products $\Pi _{n}=X_{1}X_{2}...X_{n}$ we can relax the
assumption of self-adjointness. So, suppose that $X_{i}$ are free and
identically-distributed but not necessarily self-adjoint. Also, we do not
require that $E\left( X_{i}\right) =1.$ Then we show that 
\begin{equation}
\lim_{n\rightarrow \infty }n^{-1}\log \left\Vert \Pi _{n}\right\Vert =\log 
\sqrt{E\left( X_{i}^{\ast }X_{i}\right) }.  \label{Pi_n_asymptotic}
\end{equation}

Another way to describe the behavior of $\Pi _{n}$ is to look at how the
norm of a fixed vector $\xi $ changes when we consecutively apply free
operators $X_{1},$ $...,$ $X_{n}$ to it. More precisely, suppose that the
action of the algebra of variables $X_{i}$ on a Hilbert space $H$ is
described by a cyclic representation $\pi $ and that the vector $\xi $ is
cyclic with respect to the expectation $E$. By definition, this means that $%
E\left( X\right) =\left\langle \xi ,\pi \left( X\right) \xi \right\rangle $
for every operator $X$ from a given algebra. Then we show that 
\begin{equation}
n^{-1}\log \left\Vert \pi \left( \Pi _{n}\right) \xi \right\Vert =\log \sqrt{%
E\left( X_{i}^{\ast }X_{i}\right) }.  \label{pi_Pi_asymptotics}
\end{equation}%
Note that we do not need to take the limit, since the equality holds for all 
$n.$

The reader may think of cyclic vectors as typical vectors. For example, if
the representation $\pi $ is cyclic and irreducible then cyclic vectors are
dense in $H$. In colloquial terms, (\ref{Pi_n_asymptotic}) says that for
large $n$ the product $\Pi _{n}$ cannot increase the norm of any given
vector $\xi $ by more than $\left[ E\left( X^{\ast }X\right) \right] ^{n/2}.$
And (\ref{pi_Pi_asymptotics}) says that for every cyclic vector $\xi $ this
growth rate is achieved.

One more way to capture the intuition of this result is to write 
\begin{equation*}
\lim_{n\rightarrow \infty }n^{-1}\log \left\Vert \Pi _{n}\right\Vert
=\lim_{n\rightarrow \infty }n^{-1}\log \sup_{\left\Vert x\right\Vert
=1}\left\Vert \pi \left( \Pi _{n}\right) x\right\Vert
\end{equation*}
We have shown that this limit is equal to 
\begin{equation*}
n^{-1}\log \left\Vert \pi \left( \Pi _{n}\right) \xi \right\Vert
\end{equation*}%
where $\xi $ is a cyclic vector. Thus, for large $n$ the product $\Pi _{n}$
acts uniformly in all directions. Its maximal dilation as measured by $%
\sup_{\left\Vert x\right\Vert =1}\left\Vert \pi \left( \Pi _{n}\right)
x\right\Vert $ has the same exponential order of magnitude as the dilation
in the direction of a typical vector $\xi .$

It is helpful to compare these results with the case of commutative random
variables. Suppose for the moment that $X_{i}$ are independent commutative
random variables with positive values$.$ Then, 
\begin{equation*}
\lim_{n\rightarrow \infty }n^{-1}\log \left\Vert X_{1}...X_{n}\right\Vert
=\log \left\Vert X_{i}\right\Vert ,
\end{equation*}%
where the norm of a random variable is the essential supremum norm (i.e., $%
\left\Vert X\right\Vert =\mathrm{ess}\sup_{\omega \in \Omega }\left\vert
X\left( \omega \right) \right\vert $). Indeed, for every $\varepsilon >0$
the measure of the set $\left\{ \omega :\left\vert X_{1}\left( \omega
\right) ...X_{n}\left( \omega \right) \right\vert \geq \left\Vert
X_{1}\right\Vert ...\left\Vert X_{n}\right\Vert -\varepsilon \right\} $ is
positive. Therefore $\left\Vert X_{1}...X_{n}\right\Vert =\left\Vert
X_{1}\right\Vert ^{n}.$ Note that $\log \sqrt{E\left( X_{i}^{\ast
}X_{i}\right) }\leq \log \left\Vert X_{i}\right\Vert $ and therefore the
norm of free products grows more slowly than we would expect from the
classical case.

Another interesting comparison is that with results about products of random
matrices. Let $X_{i}$ be i.i.d. random $k\times k$ matrices. Then under
suitable conditions, $\lim_{n\rightarrow \infty }n^{-1}\log \left\Vert
X_{n}...X_{1}\right\Vert $ exists almost surely. Let us denote this limit as 
$\lambda .$ \shortciteN{furstenberg63} developed a general formula for $%
\lambda ,$ and \shortciteN{cohen_newman84} derived explicit results in the
case when entries of $X_{i}$ have a joint Gaussian distribution. In
particular, if all entries of $X_{i}$ are independent and have the
distribution $\mathcal{N}\left( 0,s_{k}^{2}\right) $ then $\lambda =\left(
1/2\right) \left\{ \log \left( s_{k}^{2}\right) +\log 2+\psi \left(
k/2\right) \right\} $ where $\psi $ is the digamma function ($\psi \left(
x\right) =d\log \Gamma \left( x\right) /dx$). If the size of the matrices
grows ($k\rightarrow \infty $) then $\lambda \sim \left( 1/2\right) \log
\left( ks_{k}^{2}\right) .$ To compare this with our results, note that if $%
ks_{k}^{2}\rightarrow s^{2},$ then the sequence of random matrices
approximates a free random variable $\widetilde{X}_{i}$ with the spectral
distribution that is uniform inside the circle of radius $s.$ For this free
variable, $\ E\left( \widetilde{X}_{i}^{\ast }\widetilde{X}_{i}\right)
=s^{2},$ and our theorem shows that $\lim_{n\rightarrow \infty
}n^{-1}\left\Vert \widetilde{X}_{1}...\widetilde{X}_{n}\right\Vert =\log s.$
This limit agrees with the result for random matrices. Thus, our result can
be seen as a limiting form of results for random matrices.

The results regarding $\left\Vert Y_{n}\right\Vert $ are also interesting.
We can associate with $X_{i}$ and $Y_{n}$ probability measures $\mu _{X}$
and $\mu _{Y_{n}},$ which are called the spectral probability measures of $%
X_{i}$ and $Y_{n},$ respectively. Then the measure $\mu _{Y_{n}}$ is
determined only by $n$ and the measure $\mu _{X}$ and is called the $n$-time 
\emph{free multiplicative convolution} of $\mu _{X}$ with itself:%
\begin{equation*}
\mu _{Y_{n}}=\underset{n\text{ times}}{\underbrace{\mu _{X}\boxtimes
...\boxtimes \mu _{X}}}.
\end{equation*}%
The norm $\left\Vert Y_{n}\right\Vert $ is easy to interpret in terms of the
distribution $\mu _{Y_{n}}.$ Indeed, it is the smallest number $t$ such that
the support of $\mu _{Y_{n}}$ is inside the interval $\left[ 0,t\right] .$
Therefore, the growth in $\left\Vert Y_{n}\right\Vert $ measures the growth
in the support of the spectral probability measure if the measure is
convolved with itself using the operation of the free multiplicative
convolution.

In the case of classical multiplicative convolutions of probability
measures, the support grows exponentially, so that if $\mu _{X}$ is
supported on $\left[ 0,L_{X}\right] ,$ then the measure $\mu
_{X_{1}...X_{n}} $ is supported on $\left[ 0,\left( L_{X}\right) ^{n}\right]
.$ What we have found in the case of free multiplicative convolutions is
that if we fix $EX_{i}=1,$ then the support of the $\mu _{Y_{n}}$ grows no
faster than a linear function of $n$, i.e., the support of $\mu _{Y_{n}}$ is
inside the interval $\left[ 0,cnL_{x}\right] $ with an absolute constant $c$.

As was pointed out in the literature, a similar phenomenon occurs for sums
of free random variables. The support of measures obtained by free additive
convolutions grows much more slowly than in the case of classical additive
convolutions. This effect was called superconvergence by %
\shortciteN{bercovici_voiculescu95}. Our finding about $\left\Vert
Y_{n}\right\Vert $ can be considered as a superconvergence for free
multiplicative convolutions.

The rest of the paper is organized as follows. Section 2 formulates the
results. Section 3 contains the necessary technical background from free
probability theory. Sections 4, 5, and 6 prove the results. And Section 7
concludes.

\section{Results}

A \emph{non-commutative probability space} $\left( \mathcal{A},E\right) $ is
a unital $C^{\ast }$-algebra $\mathcal{A}$ and a positive linear functional $%
E,$ such that $E\left( I\right) =1.$ We will assume that the functional is
tracial, i.e., $E\left( AB\right) =E\left( BA\right) $ for any two operators 
$A$ and $B$ from algebra $\mathcal{A}$. The elements of algebra $\mathcal{A}$
are called \emph{random variables} and the functional $E$ is called the 
\emph{expectation}. The numbers $E\left( X^{k}\right) $ are called \emph{%
moments} of the random variable $X.$

A prototypical example of a non-commutative probability space is a group
algebra. That is, for a countable group $G$ we consider the Hilbert space $%
L^{2}\left( G,\nu \right) ,$ where $\nu $ is a counting measure, and
consider the left action of $G$ on $L^{2}\left( G,\nu \right) $: if $f\in
L^{2}\left( G,\nu \right) $ and $a,b\in G,$ then $\left[ af\right]
(b)=f(ab). $ The elements of the group algebra $\mathcal{A}$ are finite sums 
$\sum_{a\in G}x_{a}a$ and we can extend by linearity the action of the group 
$G$ on $L^{2}\left( G,\nu \right) $ to the action of the algebra $\mathcal{A}
$ on $L^{2}\left( G,\nu \right) $. We can additionally complete the
resulting operator algebra in an appropriate topology. The expectation of an
element $\sum x_{a}a$ is defined as $x_{e}$, where $e$ is the identity of
the group.

Another important example is the algebra of random matrices. The expectation
of an element $X$ in this algebra is defined as $E\left( X\right) =\mathcal{E%
}\left( N^{-1}\func{tr}\left( X\right) \right) ,$ where $\mathcal{E}$ is the
expectation with respect to underlying randomness and $N$ is the dimension
of the random matrix. For more details about these examples the reader may
consult \shortciteN{hiai_petz00}.

The concept of freeness substitutes for the concept of independence.
Consider sub-algebras. $\mathcal{A}_{1,}...,\mathcal{A}_{n}$ be given. Let $%
a_{i}$ are elements of these sub-algebras such that $a_{i}\in \mathcal{A}%
_{k\left( i\right) }.$

\begin{definition}
The algebras $\mathcal{A}_{1,}...,\mathcal{A}_{n}$ (and their elements) are 
\emph{free}, if $E\left( a_{1}...a_{m}\right) =0,$provided that $E\left(
a_{i}\right) =0$, $k(i)\neq k\left( i+1\right) $ for every $i<m,$ and $%
k\left( m\right) \neq k\left( 1\right) .$
\end{definition}

Consider the group algebra for a free group with at least two generators.
Then the operators corresponding to generators are free in the sense of the
previous definition. For the algebra of large random matrices, Voiculescu
proved the asymptotic freeness of two classically independent Gaussian
matrices, where asymptotic means that the property in the previous
definition is approached as the dimension of matrices $N\rightarrow \infty $
(see \shortciteN{voiculescu91}).

It turns out that many concepts of classical probability theory can be
transferred to the case of free random variables. For example, for a
self-adjoint variable we can define its distribution function. Indeed, if $A$
is a self-adjoint operator then by the spectral decomposition theorem it can
be written as 
\begin{equation*}
A=\int_{-\infty }^{\infty }\lambda P\left( d\lambda \right) ,
\end{equation*}%
where $P$ is a positive, projector-valued measure, i.e., a mapping that
sends sets of the real axis to orthogonal projectors$.$ This allows
definition of the \emph{spectral measure} of $A,$ $\mu _{A},$ $\ $which is a
measure with the following distribution function: 
\begin{equation*}
\mathcal{F}_{A}\left( t\right) =E\left( \int_{-\infty }^{t}P\left( d\lambda
\right) \right) .
\end{equation*}%
We can calculate the expectation of any summable function of a self-adjoint
variable $A$ by using its spectral measure:%
\begin{equation*}
Ef\left( A\right) =\int_{-\infty }^{\infty }f\left( \lambda \right) d\mu
_{A}\left( \lambda \right) .
\end{equation*}

Let $X_{1},$ $X_{2},$ ... $,$ $X_{n}$ be free identically-distributed
positive random variables. Consider $\Pi _{n}=X_{1}X_{2}...X_{n}$ and $%
Y_{n}=X_{1}\circ X_{2}\circ ...\circ X_{n}$ (by convention we multiply on
the left, so that, for example, $X_{1}\circ X_{2}\circ X_{3}\circ
X_{4}=X_{1}\circ \left( X_{2}\circ \left( X_{3}\circ X_{4}\right) \right) $
). We will see later that these variables have the same moments: $E\left(
\Pi _{n}\right) ^{k}=E\left( Y_{n}\right) ^{k}.$ As a first step let us
record some simple results about the expectation and variance of $Y_{n}$ and 
$\Pi _{n}.$ We define \emph{variance} of a random variable $A$ as%
\begin{equation*}
\sigma ^{2}\left( A\right) =:E\left( A^{\ast }A\right) -\left\vert E\left(
A\right) \right\vert ^{2}.
\end{equation*}

\begin{proposition}
\label{expectation_variance_proposition}Suppose that $X_{i}$ are
self-adjoint and $E\left( X_{i}\right) =1.$ Then $E\left( \Pi _{n}\right)
=E\left( Y_{n}\right) =1$ and $\sigma ^{2}\left( \Pi _{n}\right) =\sigma
^{2}\left( Y_{n}\right) =n\sigma ^{2}\left( X_{i}\right) .$
\end{proposition}

Note that the linear growth in the variance of $\Pi _{n}=X_{1}...X_{n}$ is
in contrast with the classical case, where only the variance of $\log \left(
X_{1}...X_{n}\right) $ grows linearly. We will prove this Proposition later
when we have more technical tools available. Before that we are going to
formulate the main results.

Let $\left\Vert A\right\Vert $ denote the usual operator norm of operator $A$%
.

\begin{theorem}
\label{theorem1}Suppose that $X_{1}$, ..., $X_{n}$ \ are
identically-distributed positive self-adjoint free variables$.$ Suppose also
that $E\left( X_{i}\right) =1.$ Then\newline
(1) there exists such a constant, $c,$ that $\left\Vert Y_{n}\right\Vert
\leq c\left\Vert X_{i}\right\Vert n;$ \newline
and\newline
(2) $\left\Vert Y_{n}\right\Vert \geq \sigma \left( X_{i}\right) \sqrt{n}.$
\end{theorem}

For the next theorem define 
\begin{equation*}
\gamma =\sigma \left( \frac{X_{i}^{\ast }X_{i}}{E\left( X_{i}^{\ast
}X_{i}\right) }\right) \geq 0
\end{equation*}

\begin{theorem}
\label{theorem2} Suppose that $X_{1}$, ..., $X_{n}$ \ are free
identically-distributed variables (not necessarily self-adjoint)$.$ Then 
\newline
(1) there exists such a constant, $c,$ that $\left\Vert \Pi _{n}\right\Vert
\leq c\left\Vert X_{i}\right\Vert \sqrt{n}\left[ E\left( X_{i}^{\ast
}X_{i}\right) \right] ^{\left( n-1\right) /2};$\newline
and\newline
(2) $\left\Vert \Pi _{n}\right\Vert \geq \gamma ^{1/2}n^{1/4}\left[ E\left(
X_{i}^{\ast }X_{i}\right) \right] ^{n/2}.$
\end{theorem}

\begin{corollary}
\label{corollary1}Suppose that $X_{1}$, ..., $X_{n}$ \ are free
identically-distributed variables (not necessarily self-adjoint)$.$ Then 
\begin{equation*}
\lim_{n\rightarrow \infty }n^{-1}\log \left\Vert \Pi _{n}\right\Vert =\log 
\sqrt{E\left( X_{i}^{\ast }X_{i}\right) }
\end{equation*}
\end{corollary}

Next, suppose that the algebra $\mathcal{A}$ acts on an
(infinitely-dimensional) Hilbert space $H.$ In other words, let $\pi $ be a
representation of $\mathcal{A}$. We call representation $\pi $ \emph{cyclic }%
if there exists such a vector $\xi \in H$ that $E\left( X\right)
=\left\langle \xi ,\pi \left( X\right) \xi \right\rangle $ for all operators 
$X\in \mathcal{A}.$ The vectors with this property are also called \emph{%
cyclic}.

\begin{theorem}
\label{theorem3}Suppose $\pi $ is a cyclic representation of $\mathcal{A}$, $%
\xi $ is its cyclic vector, and $X_{1}$, ..., $X_{n}$ \ are free
identically-distributed variables from $\mathcal{A}$. Then 
\begin{equation*}
n^{-1}\log \left\Vert \pi \left( \Pi _{n}\right) \xi \right\Vert =\log \sqrt{%
E\left( X_{i}^{\ast }X_{i}\right) }
\end{equation*}
\end{theorem}

\begin{corollary}
If $\pi $ and $\xi $ are cyclic then 
\begin{equation*}
\log \left\Vert \Pi _{n}\right\Vert \sim \log \left\Vert \pi \left( \Pi
_{n}\right) \xi \right\Vert \sim n\log \left\Vert \pi \left( X_{1}\right)
\xi \right\Vert
\end{equation*}
as $n\rightarrow \infty .$
\end{corollary}

\section{Preliminaries}

The \emph{Cauchy transform} of a bounded random variable $A$ is defined as
follows:%
\begin{equation*}
G_{A}\left( z\right) =E\left( \frac{1}{z-A}\right) =\frac{1}{z}%
+\sum_{k=1}^{\infty }\frac{E\left( A^{k}\right) }{z^{k+1}}.
\end{equation*}%
This power series is convergent for $\left\vert z\right\vert >\left\Vert
A\right\Vert .$ Let us also define the $\psi $\emph{-function} of $A$:%
\begin{equation*}
\psi _{A}\left( z\right) =E\left( \frac{1}{1-zA}\right)
-1=\sum_{k=1}^{\infty }E\left( A^{k}\right) z^{k}.
\end{equation*}%
The $\psi $-function is convergent for $\left\vert z\right\vert \leq
\left\Vert A\right\Vert ^{-1}$ and it is related to the Cauchy transform by
the following equality: 
\begin{equation*}
G_{n}\left( z\right) =z^{-1}\left[ \psi _{n}\left( z^{-1}\right) +1\right] .
\end{equation*}

If $A$ is bounded and $E\left( A\right) \neq 0,$ then for $z$ in a
sufficiently small neighborhood of $0,$ the inverse of $\psi _{A}\left(
z\right) $ is defined, which we denote as $\psi _{A}^{-1}\left( z\right) .$
Then the \emph{S-transform} is defined as 
\begin{equation}
S_{A}\left( z\right) =\left( 1+\frac{1}{z}\right) \psi _{A}^{-1}\left(
z\right)  \label{S_definition}
\end{equation}

Let us write out several first terms in the power expansions for $\psi
\left( z\right) ,$ $\psi ^{-1}\left( z\right) ,$ and $S\left( z\right) .$
Suppose for simplicity that $E\left( A\right) =1$ and let $E\left(
A^{k}\right) =m_{k}$. Then,%
\begin{eqnarray*}
\psi \left( z\right) &=&z+m_{2}z^{2}+m_{3}z^{3}+..., \\
\psi ^{-1}\left( z\right) &=&z-m_{2}z^{2}-\left( m_{3}-2m_{2}^{2}\right)
z^{3}+..., \\
S\left( z\right) &=&1+(1-m_{2})z+\left( 2m_{2}^{2}-m_{2}-m_{3}\right)
z^{2}+...
\end{eqnarray*}

The main theorem regarding the multiplication of free random variables was
proved by \shortciteN{voiculescu87}. Later the proof was significantly
simplified by \shortciteN{haagerup97}.

\begin{theorem}[Voiculescu]
\label{Voiculescu_multiplication}Suppose $X$ and $Y$ are bounded free random
variables. Suppose also that $E\left( X\right) \neq 0$ and $E\left( Y\right)
\neq 0.$ Then 
\begin{equation*}
S_{XY}\left( z\right) =S_{X}\left( z\right) S_{Y}\left( z\right) .
\end{equation*}
\end{theorem}

In particular, this theorem implies that $S_{\Pi _{n}}=S_{Y_{n}}=\left(
S_{X}\right) ^{n},$ where $S_{X}$ denotes the $S$-transform of any of $%
X_{i}. $ Now it is easy to prove Proposition \ref%
{expectation_variance_proposition}. Indeed, let us denote $S_{\Pi _{n}}$ as $%
S_{n}.$ Then, using the power expansions we can write:%
\begin{eqnarray*}
S_{n}\left( z\right) &=&1+\left( 1-m_{2}^{(n)}\right) z+... \\
&=&\left( S_{X}\right) ^{n}=1+n\left( 1-m_{2}\right) z+...,
\end{eqnarray*}%
where $m_{2}^{\left( n\right) }=:E\left( \Pi _{n}\right) ^{2}$ and $%
m_{2}=:E\left( X_{i}\right) ^{2}.$ Then, using power expansion in (\ref%
{S_definition}), we conclude that $E\left( \Pi _{n}\right) =1.$ Next, by
definition, $\sigma ^{2}\left( X_{i}\right) =m_{2}-1$ and $\sigma ^{2}\left(
\Pi _{n}\right) =m_{2}^{\left( n\right) }-1.$ Therefore, we can conclude
that $\sigma ^{2}\left( \Pi _{n}\right) =n\sigma ^{2}\left( X\right) .$ QED.

\section{Proof of Theorem \protect\ref{theorem1}}

Throughout this section we assume that $X_{i}$ are self-adjoint, $E\left(
X_{i}\right) =1,$ and the support of the spectral distribution of $X_{i}$
belongs to $\left[ 0,L\right] .$

Let us first go in a simpler direction and derive a lower bound on $%
\left\Vert Y_{n}\right\Vert .$ That is, we are going to prove claim (2) of
the theorem. From Proposition \ref{expectation_variance_proposition}, we
know that $E\left( Y_{n}\right) =1$ and $\sigma ^{2}\left( Y_{n}\right)
=n\sigma ^{2}\left( X_{i}\right) .$ It is clear that for every positive
random variable $A,$ it is true that $E\left( A^{2}\right) \leq \left\Vert
A\right\Vert ^{2}$ and therefore $\left\Vert A\right\Vert \geq \sqrt{\sigma
^{2}\left( A\right) +\left[ E\left( A\right) \right] ^{2}}.$ Applying this
to $Y_{n}$, we get $\left\Vert Y_{n}\right\Vert \geq \sqrt{n\sigma ^{2}+1}.$
In particular, $\left\Vert Y_{n}\right\Vert >\sigma \sqrt{n},$ so (2) is
proved.

Now let us prove claim (1). By Theorem \ref{Voiculescu_multiplication}, $%
S_{n}\left( z\right) =\left( S_{X}\left( z\right) \right) ^{n}$. The idea of
the proof is to investigate how $\left\vert S_{X}\left( z\right) \right\vert
^{n}$ behaves for small $z.$ It turns out that if $z$ is of the order of $%
n^{-1}$, then $\left\vert S_{X}\left( z\right) \right\vert ^{n}>c$ where $c$
is a constant that does not depend on $n.$ We will show that this fact
implies that $\psi _{n}\left( z\right) $ (i.e., the $\psi $-function for $%
Y_{n}$) has the convergent power series in the area $\left\vert z\right\vert
<\left( cn\right) ^{-1}$ and that therefore the Cauchy transform of $Y_{n}$
has the convergent power series in $\left\vert z\right\vert >cn.$ This fact
and the Perron-Stieltjes inversion formula imply that the support of the
distribution of $Y_{n}$ is inside $\left[ -cn,cn\right] .$

In the proof we need the result about functional inversions formulated
below. By a function holomorphic in a domain, $D,$ we mean a function which
is bounded and differentiable in $D.$

\begin{lemma}[Lagrange's inversion formula]
\label{Lagranges_series_around0} Suppose $f$ is a function of a complex
variable, which is holomorphic in a neighborhood of $z_{0}=0$ and has the
Taylor expansion 
\begin{equation*}
f(z)=z+\sum_{k=2}^{\infty }a_{k}z^{k}+...,
\end{equation*}%
converging for all sufficiently small $z.$ Then the functional inverse of $%
f\left( z\right) $ is well defined in a neighborhood of $0$ and the Taylor
series of the inverse is given by the following formula:%
\begin{equation*}
f^{-1}\left( u\right) =u+\sum_{k=2}^{\infty }\left[ \frac{1}{2\pi ik}%
\oint_{\gamma }\frac{dz}{f(z)^{k}}\right] u^{k},
\end{equation*}%
where $\gamma $ is a circle around $0,$ in which $f$ has only one zero..
\end{lemma}

For the proof see Theorems II.3.2 and II.3.3 in \shortciteN{markushevich77},
or Section 7.32 in \shortciteN{whittaker_watson27}.

\begin{lemma}
\label{moment_bound} $E\left( X^{k}\right) \leq L^{k-1}.$
\end{lemma}

\textbf{Proof:} 
\begin{equation*}
E\left( X^{k}\right) =\int_{0}^{L}\lambda ^{k}d\mu _{X}\left( \lambda
\right) \leq L^{k-1}\int_{0}^{L}\lambda d\mu _{X}\left( \lambda \right)
=L^{k-1},
\end{equation*}%
where $d\mu _{X}$ denotes the spectral distribution of the variable $X.$ QED.

\begin{lemma}
The function $\psi _{X}\left( z\right) $ is has only one zero in $\left\vert
z\right\vert \leq \left( 4L\right) ^{-1},$ and if $\left\vert z\right\vert
=\left( 4L\right) ^{-1},$ then $\left\vert \psi _{X}\left( z\right)
\right\vert \geq \left( 6L\right) ^{-1}.$
\end{lemma}

\textbf{Proof:} If $\left\vert z\right\vert \leq \left( 4L\right) ^{-1}$
then 
\begin{eqnarray*}
\left\vert \psi _{X}\left( z\right) -z\right\vert &\leq &\left\vert
z\right\vert \sum_{k=2}^{\infty }E\left( X^{k}\right) \left\vert
z\right\vert ^{k-1} \\
&\leq &\left\vert z\right\vert \sum_{k=1}^{\infty }\frac{1}{4^{k}}=\frac{%
\left\vert z\right\vert }{3}.
\end{eqnarray*}%
Therefore, by Rouch\'{e}'s theorem, $\psi _{X}\left( z\right) $ has only one
zero in this area.

If $\left\vert z\right\vert =\left( 4L\right) ^{-1},$ then 
\begin{eqnarray*}
\left\vert \psi _{X}\left( z\right) \right\vert &\geq &\left\vert
z\right\vert -\sum_{k=2}^{\infty }E\left( X^{k}\right) \left\vert
z\right\vert ^{k} \\
&\geq &\left\vert z\right\vert \left( 1-\sum_{k=1}^{\infty }\frac{1}{4^{k}}%
\right) \\
&=&\frac{1}{4L}\left( 1-\frac{1}{3}\right) =\frac{1}{6L}.
\end{eqnarray*}

QED.

By Lemma \ref{Lagranges_series_around0}, we can expand the functional
inverse of $\psi _{X}\left( z\right) $ as follows:%
\begin{equation*}
\psi _{X}^{-1}\left( u\right) =u+\sum_{k=2}^{\infty }c_{k}u^{k},
\end{equation*}%
where 
\begin{equation*}
c_{k}=\frac{1}{2\pi ik}\int_{\gamma }\frac{dz}{\left[ \psi _{X}\left(
z\right) \right] ^{k}}
\end{equation*}

\begin{lemma}
If $\left\vert u\right\vert \leq \left( 72Ln\right) ^{-1},$ then 
\begin{equation*}
\left\vert \frac{\psi _{X}^{-1}\left( u\right) }{u}-1\right\vert \leq \frac{1%
}{7n}.
\end{equation*}
\end{lemma}

\textbf{Proof: }Using the previous lemma we can estimate $c_{k}$:%
\begin{equation*}
c_{k}\leq \frac{1}{k}\frac{1}{4L}\left( 6L\right) ^{k}\leq \frac{3}{2}\left(
6L\right) ^{k-1}.
\end{equation*}%
Then 
\begin{eqnarray*}
\left\vert \frac{\psi _{X}^{-1}\left( u\right) }{u}-1\right\vert
&=&\left\vert \sum_{k=2}^{\infty }c_{k}u^{k-1}\right\vert \\
&\leq &\frac{3}{2}\dsum\limits_{k=1}^{\infty }\left( \frac{1}{12n}\right)
^{k}=\frac{3}{2}\frac{1}{12n-1} \\
&=&\frac{3}{2}\frac{12n}{12n-1}\frac{1}{12n}\leq \frac{1}{7n},
\end{eqnarray*}%
provided that $\left\vert u\right\vert \leq \left( 72Ln\right) ^{-1}.$ QED.

\begin{lemma}
If $\left\vert u\right\vert \leq \left( 72Ln\right) ^{-1},$ then%
\begin{equation*}
\left\vert 1-S_{X}\left( u\right) \right\vert \leq \frac{1}{6n}.
\end{equation*}
\end{lemma}

\textbf{Proof:} Recall that $S_{X}\left( u\right) =\left( 1+u\right) \psi
_{X}^{-1}\left( u\right) /u.$ Then we can write:%
\begin{eqnarray*}
\left\vert 1-S_{X}\left( u\right) \right\vert &=&\left\vert u+\left(
1+u\right) \left( \frac{\psi _{X}^{-1}\left( u\right) }{u}-1\right)
\right\vert \\
&\leq &\left\vert u\right\vert +\left\vert 1+u\right\vert \left\vert \frac{%
\psi _{X}^{-1}\left( u\right) }{u}-1\right\vert .
\end{eqnarray*}%
Then the previous lemma implies that for $\left\vert u\right\vert \leq
\left( 72Ln\right) ^{-1}$ and $n\geq 2,$ we have the estimate: 
\begin{equation*}
\left\vert 1-S_{X}\left( u\right) \right\vert \leq \frac{1}{72Ln}+\left\vert
1+\frac{1}{72Ln}\right\vert \frac{1}{7n}.
\end{equation*}%
Note that $L\geq 1$ because $EX=1.$ Therefore, 
\begin{equation*}
\left\vert 1-S_{X}\left( u\right) \right\vert \leq \frac{1}{72n}+\frac{73}{72%
}\frac{1}{7n}\leq \frac{1}{6n}.
\end{equation*}%
QED.

\begin{lemma}
\label{estimate_Sn_lemma} For all positive integer $n$ if $\left\vert
u\right\vert \leq \left( 72Ln\right) ^{-1},$ then 
\begin{equation*}
e^{1/6}\geq \left\vert S_{X}\left( u\right) \right\vert ^{n}\geq e^{-1/3}.
\end{equation*}
\end{lemma}

\textbf{Proof:} Let us first prove the upper bound on $\left\vert
S_{X}\left( u\right) \right\vert ^{n}.$ The previous lemma implies that 
\begin{equation*}
\left\vert S_{X}\left( u\right) \right\vert ^{n}\leq \left( 1+\frac{1}{6n}%
\right) ^{n}\leq e^{1/6}.
\end{equation*}%
Now let us prove the lower bound. The previous lemma implies that 
\begin{equation*}
\left\vert S_{X}\left( u\right) \right\vert ^{n}\geq \left( 1-\frac{1}{6n}%
\right) ^{n}.
\end{equation*}%
In an equivalent form, 
\begin{equation}
n\log \left\vert S_{X}\left( u\right) \right\vert \geq n\log \left( 1-\frac{1%
}{6n}\right) .  \label{estimate_Sn}
\end{equation}%
Recall the following elementary inequality: If $x\in \left[ 0,1-e^{-1}\right]
,$ then 
\begin{equation*}
\log \left( 1-x\right) \geq -2x.
\end{equation*}%
Let $x=1/\left( 6n\right) $. Then 
\begin{equation*}
\log \left( 1-\frac{1}{6n}\right) \geq -\frac{1}{3n}.
\end{equation*}%
Substituting this in (\ref{estimate_Sn}), we get 
\begin{equation*}
n\log \left\vert S_{X}\left( u\right) \right\vert \geq -\frac{1}{3},
\end{equation*}%
or 
\begin{equation*}
\left\vert S_{X}\left( u\right) \right\vert ^{n}\geq e^{-1/3}.
\end{equation*}%
QED.

By Theorem $\ $\ref{Voiculescu_multiplication},\ $S_{n}(u)=:\left[
S_{X}\left( u\right) \right] ^{n}$ is the $S$-transform of the variable $%
Y_{n}.$ The corresponding inverse $\psi $-function is $\psi _{n}^{-1}\left(
u\right) =uS_{n}\left( u\right) /\left( 1+u\right) .$

First, we estimate $S_{n}\left( u\right) -1$.

\begin{lemma}
\label{lemma_Sn_estimate}If $\left\vert u\right\vert \leq \left( 72Ln\right)
^{-1},$ then%
\begin{equation*}
\left\vert S_{n}\left( u\right) -1\right\vert \leq \frac{1}{5}.
\end{equation*}
\end{lemma}

\textbf{Proof: }Write 
\begin{eqnarray*}
\left\vert S_{X}\left( u\right) ^{n}-1\right\vert &\leq &\left\vert
S_{X}\left( u\right) -1\right\vert \left( \left\vert S_{X}\left( u\right)
\right\vert ^{n-1}+\left\vert S_{X}\left( u\right) \right\vert
^{n-2}+...+1\right) \\
&\leq &\frac{1}{6n}e^{1/6}n\leq \frac{1}{5}.
\end{eqnarray*}

QED.

\begin{lemma}
\label{estimate_psi_n_lemma} The function $\psi _{n}^{-1}\left( u\right) $
has only one zero in $\left\vert u\right\vert =\left( 72Ln\right) ^{-1}$ and
if $\left\vert u\right\vert =\left( 72Ln\right) ^{-1},$ then%
\begin{equation*}
\left\vert \psi _{n}^{-1}\left( u\right) \right\vert \geq \frac{1}{102Ln}.
\end{equation*}
\end{lemma}

\textbf{Proof: }Recall that by definition in (\ref{S_definition}), $\psi
_{n}^{-1}\left( u\right) =uS_{n}\left( u\right) /\left( 1+u\right) .$
Therefore, 
\begin{equation*}
\left\vert \psi _{n}^{-1}\left( u\right) -u\right\vert =\left\vert
u\right\vert \left\vert \frac{S_{n}\left( u\right) -\left( 1+u\right) }{1+u}%
\right\vert
\end{equation*}%
and by Lemma \ref{lemma_Sn_estimate} we have the following estimate:%
\begin{eqnarray*}
\left\vert \frac{S_{n}\left( u\right) -\left( 1+u\right) }{1+u}\right\vert
&\leq &\frac{1}{1-\left\vert u\right\vert }\left\vert S_{n}\left( u\right)
-1\right\vert +\frac{\left\vert u\right\vert }{1-\left\vert u\right\vert } \\
&\leq &\frac{72}{71}\frac{1}{5}+\frac{1}{71}\leq \frac{1}{4}.
\end{eqnarray*}%
Therefore, by Rouch\'{e}'s theorem, $\psi _{n}^{-1}\left( u\right) $ has
only one zero in $\left\vert u\right\vert \leq \left( 72Ln\right) ^{-1}.$

Next, note that $\psi _{n}^{-1}\left( u\right) =uS_{n}\left( u\right)
/\left( 1+u\right) $ and if $\left\vert u\right\vert =\left( 72Ln\right)
^{-1},$ then%
\begin{equation*}
\left\vert \frac{u}{1+u}\right\vert \geq \frac{1}{72Ln}/\left( 1+\frac{1}{%
72Ln}\right) \geq \frac{1}{73Ln}.
\end{equation*}%
Using Lemma \ref{estimate_Sn_lemma}, we get:%
\begin{equation*}
\left\vert \psi _{n}^{-1}\left( u\right) \right\vert \geq \frac{1}{73Ln}%
e^{-1/3}\geq \frac{1}{102Ln}.
\end{equation*}%
QED.

Now we again apply Lemma \ref{Lagranges_series_around0} and obtain the
following formula: 
\begin{equation}
\psi _{n}\left( z\right) =z+\sum_{k=2}^{\infty }\left[ \frac{1}{2\pi ik}%
\oint_{\gamma }\frac{du}{\left[ \psi _{n}^{-1}\left( u\right) \right] ^{k}}%
\right] z^{k},  \label{psi_series}
\end{equation}%
where we can take the circle $\left\vert u\right\vert =\left( 72Ln\right)
^{-1}$ as $\gamma .$

\begin{lemma}
\label{psi_convergence_radius}The radius of convergence of series (\ref%
{psi_series}) is at least $\left( 102Ln\right) ^{-1}.$
\end{lemma}

\textbf{Proof:} By the previous lemma, the coefficient before $z^{k}$ can be
estimated as follows:%
\begin{equation*}
\left\vert c_{k}\right\vert \leq \frac{1}{k}\frac{1}{72Ln}\left(
102Ln\right) ^{k}.
\end{equation*}%
This implies that series (\ref{psi_series}) converges at least for $%
\left\vert z\right\vert \leq \left( 102Ln\right) ^{-1}.$ QED.

\begin{lemma}
\label{distribution_support}The support of the spectral distribution of $%
Y_{n}=X_{1}\circ X_{2}\circ ...\circ X_{n}$ belongs to the interval $\left[
-102Ln,102Ln\right] .$
\end{lemma}

\textbf{Proof:} The variable $Y_{n}$ is self-adjoint \ and has a
well-defined spectral measure, $\mu _{n}\left( dx\right) ,$ supported on the
real axis. We can infer the Cauchy transform of this measure from $\psi
_{n}\left( z\right) $:%
\begin{equation*}
G_{n}\left( z\right) =z^{-1}\left[ \psi _{n}\left( z^{-1}\right) +1\right] .
\end{equation*}%
Using Lemma \ref{psi_convergence_radius}, we can conclude that the power
series for $G_{n}\left( z\right) $ around $z=\infty $ converges in the area $%
\left\vert z\right\vert >102Ln.$ The coefficients of this series are real.
Therefore, using the Perron-Stieltjes formula we conclude that $\mu
_{n}\left( dx\right) $ is zero outside of the interval $\left[ -102Ln,102Ln%
\right] .$ QED.

Lemma \ref{distribution_support} implies the statement of Theorem \ref%
{theorem1}.

\section{Proof of Theorem \protect\ref{theorem2}}

The norm of the operator $\Pi _{n}$ coincides with the square root of the
norm of the operator $\Pi _{n}^{\ast }\Pi _{n}.$ Therefore, all we need to
do is to estimate the norm of the self-adjoint operator $\Pi _{n}^{\ast }\Pi
_{n}.$

\begin{lemma}
\label{XstarX}For every bounded operator $X\in \mathcal{A},$ products $%
X^{\ast }X$ and $XX^{\ast }$ have the same spectral distribution.
\end{lemma}

\textbf{Proof:} Since $E$ is tracial, $E\left( X^{\ast }X\right)
^{k}=E\left( XX^{\ast }\right) ^{k}$. Therefore, $X^{\ast }X$ and $XX^{\ast
} $ have the same sequence of moments and, therefore, the same distribution.
QED.

If two variables $A$ and $B$ have the same sequence of moments, we say that
they are \emph{equivalent} and write $A\sim B.$ In particular, two
self-adjoint bounded variables have the same spectral distribution if and
only if they are equivalent.

\begin{lemma}
\label{equivalence_of_products} Let $A,$ $B,$ and $C$ be three bounded
operators from a non-commutative probability space $\mathcal{A}$. If $A\sim
B $, $A$ is free from $C$, and $B$ is free from $C,$ then $A+C\sim B+C$, $%
AC\sim BC$, $\ $and $CA\sim CB.$
\end{lemma}

\textbf{Proof: }Since $A$ and $C$ are free, the moments of $A+C$ can be
computed from the moments of $A$ and $C.$ The computation is exactly the
same for $B+C,$ since $B$ and $C$ are also free. In addition, we know that $%
A $ and $B$ have the same moments. Consequently, $A+C$ has the same moments
as $B+C,$ i.e., $A+C\sim B+C$. The other equivalences are obtained
similarly. QED.

\begin{lemma}
\label{equivalence_and_Stransform}If $A\sim B$, then $S_{A}\left( z\right)
=S_{B}\left( z\right) $. In words, if two variables are equivalent, then
they have the same $S$-transform.
\end{lemma}

\textbf{Proof:} From the definition of the $\psi $-function, it is clear
that if $A\sim B$, then $\psi _{A}\left( z\right) =\psi _{B}\left( z\right) $%
. This implies that $\psi _{A}^{-1}\left( z\right) =\psi _{B}^{-1}\left(
z\right) $ and therefore $S_{A}\left( z\right) =S_{B}\left( z\right) .$ QED.

For example, $S_{X_{i}^{\ast }X_{i}}\left( z\right) $ does not depend on $%
i\,\ $\ and we will denote this function as $S_{X^{\ast }X}\left( z\right) .$

\begin{lemma}
\label{lemma_Pi_n_star_Pi_n}If $X_{1},$ $...,$ $X_{n}$ are free, then 
\begin{equation*}
\Pi _{n}^{\ast }\Pi _{n}\sim X_{n}^{\ast }X_{n}...X_{1}^{\ast }X_{1}
\end{equation*}%
and if $X_{1},$ $...,$ $X_{n}$ are in addition identically distributed, then 
\begin{equation*}
S_{\Pi _{n}^{\ast }\Pi _{n}}=S_{\Pi _{n}\Pi _{n}^{\ast }}=\left( S_{X^{\ast
}X}\right) ^{n}
\end{equation*}
\end{lemma}

\textbf{Proof: }We will use induction. For $n=1,$ we have $\Pi _{1}^{\ast
}\Pi _{1}=X_{1}^{\ast }X_{1}$. Therefore $S_{\Pi _{1}^{\ast }\Pi
_{1}}=S_{X^{\ast }X}.$ Suppose that the statement is proved for $n-1.$ Then 
\begin{eqnarray*}
\Pi _{n}^{\ast }\Pi _{n} &=&X_{n}^{\ast }...X_{1}^{\ast }X_{1}...X_{n} \\
&\sim &X_{n}X_{n}^{\ast }X_{n-1}^{\ast }...X_{1}^{\ast }X_{1}...X_{n-1},
\end{eqnarray*}%
where the equivalence holds because $E$ is tracial and it is easy to check
that the products have the same moments. Therefore, 
\begin{eqnarray*}
\Pi _{n}^{\ast }\Pi _{n} &\sim &\left( X_{n}X_{n}^{\ast }\right) \Pi
_{n-1}^{\ast }\Pi _{n-1} \\
&\sim &\left( X_{n}^{\ast }X_{n}\right) \Pi _{n-1}^{\ast }\Pi _{n-1}
\end{eqnarray*}%
by Lemmas \ref{XstarX} and \ref{equivalence_of_products}. Then the inductive
hypothesis implies that 
\begin{equation*}
\Pi _{n}^{\ast }\Pi _{n}\sim X_{n}^{\ast }X_{n}...X_{1}^{\ast }X_{1}.
\end{equation*}

Using Lemma \ref{equivalence_and_Stransform} and Theorem \ref%
{Voiculescu_multiplication}, we write: 
\begin{equation*}
S_{\Pi _{n}^{\ast }\Pi _{n}}=\left( S_{X^{\ast }X}\right) ^{n}.
\end{equation*}%
Since $\Pi _{n}^{\ast }\Pi _{n}\sim \Pi _{n}\Pi _{n}^{\ast },$ therefore, $%
S_{\Pi _{n}^{\ast }\Pi _{n}}=S_{\Pi _{n}\Pi _{n}^{\ast }}=\left( S_{X^{\ast
}X}\right) ^{n}.$ QED.

We have managed to represent $S_{\Pi _{n}^{\ast }\Pi _{n}}$ as $\left(
S_{X^{\ast }X}\right) ^{n}$ and therefore all the arguments of the previous
section are applicable, except that we are interested in $\left( S_{X^{\ast
}X}\right) ^{n}$ rather than in $\left( S_{X}\right) ^{n}.$ In particular,
we can conclude that the following lemma holds:

\begin{lemma}
\label{PI_PI_support_lemma} Define%
\begin{equation*}
\gamma =\sigma \left( \frac{X_{i}^{\ast }X_{i}}{E\left( X_{i}^{\ast
}X_{i}\right) }\right) .
\end{equation*}%
Then\newline
(1) $\left\Vert \Pi _{n}^{\ast }\Pi _{n}\right\Vert \leq 102\left\Vert
X_{i}\right\Vert ^{2}nE\left( X_{i}^{\ast }X_{i}\right) ^{n-1},$ and \newline
(2) $\left\Vert \Pi _{n}^{\ast }\Pi _{n}\right\Vert \geq \gamma \sqrt{n}%
E\left( X_{i}^{\ast }X_{i}\right) ^{n}.$
\end{lemma}

\textbf{Proof: }Let us introduce variables $R_{i}=s^{-1}X_{i},$ where $%
s^{2}=E\left( X^{\ast }X\right) .$ Then $\left\Vert R_{i}^{\ast
}R_{i}\right\Vert =\left( \left\Vert X_{i}\right\Vert /s\right) ^{2}$ and $%
E\left( R_{i}^{\ast }R_{i}\right) =1.$ Let $\widetilde{\Pi }%
_{n}=R_{1}...R_{n}.$ Then $\Pi _{n}^{\ast }\Pi _{n}=s^{2n}\widetilde{\Pi }%
_{n}^{\ast }\widetilde{\Pi }_{n}$ and the $S$-transform of $\widetilde{\Pi }%
_{n}^{\ast }\widetilde{\Pi }_{n}$ is $\left( S_{R^{\ast }R}\right) ^{n}.$

Note that $\widetilde{\Pi }_{n}^{\ast }\widetilde{\Pi }_{n}$ has the same $S$%
-transform and therefore the same spectral distribution as $\left(
R_{1}^{\ast }R_{1}\right) \circ ...\circ \left( R_{n}^{\ast }R_{n}\right) .$
Using Theorem \ref{theorem1}, we conclude that $\left\Vert \widetilde{\Pi }%
_{n}^{\ast }\widetilde{\Pi }_{n}\right\Vert \leq 102\left( \left\Vert
X_{i}\right\Vert /s\right) ^{2}n.$ It follows that $\left\Vert \Pi
_{n}^{\ast }\Pi _{n}\right\Vert \leq 102\left\Vert X_{i}\right\Vert
^{2}s^{2n-2}n.$ In addition, Theorem \ref{theorem1} implies that 
\begin{eqnarray*}
\left\Vert \widetilde{\Pi }_{n}^{\ast }\widetilde{\Pi }_{n}\right\Vert &>&%
\sqrt{n}\sigma \left( R_{i}^{\ast }R_{i}\right) \\
&=&\gamma \sqrt{n}
\end{eqnarray*}%
Consequently, 
\begin{equation*}
\left\Vert \Pi _{n}^{\ast }\Pi _{n}\right\Vert \geq \gamma \sqrt{n}s^{2n}.
\end{equation*}%
QED.

From Lemma \ref{PI_PI_support_lemma} we conclude that%
\begin{equation*}
\left\Vert \Pi _{n}\right\Vert \leq 11\left\Vert X_{i}\right\Vert \sqrt{n}%
\left[ E\left( X_{i}^{\ast }X_{i}\right) \right] ^{\left( n-1\right) /2},
\end{equation*}%
and 
\begin{equation*}
\left\Vert \Pi _{n}\right\Vert \geq \gamma ^{1/2}n^{1/4}\left[ E\left(
X_{i}^{\ast }X_{i}\right) \right] ^{n/2}
\end{equation*}%
This completes the proof of Theorem \ref{theorem2}.

\section{Proof of Theorem \protect\ref{theorem3}}

By definition of the cyclic vector, we have:%
\begin{eqnarray*}
\left\Vert \pi \left( \Pi _{n}\right) \xi \right\Vert ^{2} &=&\left\langle
\pi \left( \Pi _{n}\right) \xi ,\pi \left( \Pi _{n}\right) \xi \right\rangle
\\
&=&\left\langle \xi ,\pi \left( \Pi _{n}^{\ast }\Pi _{n}\right) \xi
\right\rangle \\
&=&E\left( \Pi _{n}^{\ast }\Pi _{n}\right) .
\end{eqnarray*}

Using Lemma \ref{lemma_Pi_n_star_Pi_n}, we continue this as follows:%
\begin{eqnarray*}
E\left( \Pi _{n}^{\ast }\Pi _{n}\right) &=&E\left( X_{n}^{\ast
}X_{n}...X_{1}^{\ast }X_{1}\right) \\
&=&\left[ E\left( X^{\ast }X\right) \right] ^{n}.
\end{eqnarray*}%
Consequently, 
\begin{equation*}
n^{-1}\log \left\Vert \Pi _{n}\xi \right\Vert =\frac{1}{2}\log E\left(
X^{\ast }X\right) .
\end{equation*}%
QED.

\section{Concluding Remarks}

We have investigated how the norms of $\Pi _{n}=X_{1}...X_{n}$ and $%
Y_{n}=X_{1}\circ ...\circ X_{n}$ grow as $n\rightarrow \infty .$ For $%
\left\Vert \Pi _{n}\right\Vert ,$ we have shown that $\lim_{n\rightarrow
\infty }n^{-1}\log \left\Vert \Pi _{n}\right\Vert $ exists and equals $\log 
\sqrt{E\left( X_{i}^{\ast }X_{i}\right) }.$ For $\left\Vert Y_{n}\right\Vert
,$ we have proved that the growth rate of $\left\Vert Y_{n}\right\Vert $ is
somewhere between $\sqrt{n}$ and $n.$ There remains the question of whether $%
\lim_{n\rightarrow \infty }n^{-s}\left\Vert Y_{n}\right\Vert $ exists for
some $s.$

Another interesting question, which is not resolved in this paper, is how
the spectral radius of $\Pi _{n}$ grows. \ Indeed, for $Y_{n}$, the norm
coincides with the spectral radius. But for $\Pi _{n},$ the norm and the
spectral radius are different because $\Pi _{n}$ is not self-adjoint.

\bibliographystyle{CHICAGO}
\bibliography{comtest}

\end{document}